\documentclass[reqno]{amsart}
\usepackage[left=1in,right=1in,top=.95in,bottom=.95in]{geometry}
\usepackage{tikz}
\usetikzlibrary{shapes,snakes,calc,arrows}
\usepackage{amsthm}
\usepackage{amscd}
\usepackage{amsfonts}
\usepackage{amsmath}
\usepackage{amssymb}
\usepackage{amsrefs}
\usepackage{mathrsfs}
\usepackage{multirow}
\usepackage{verbatim}
\usepackage{url}
\usepackage{graphicx}

\begin{document}


\newcommand{\supp}{\text{supp}}
\newcommand{\Aut}{\text{Aut}}
\newcommand{\Gal}{\text{Gal}}
\newcommand{\Inn}{\text{Inn}}
\newcommand{\Irr}{\text{Irr}}
\newcommand{\Ker}{\text{Ker}}
\newcommand{\N}{\mathbb{N}}
\newcommand{\Z}{\mathbb{Z}}
\newcommand{\Q}{\mathbb{Q}}
\newcommand{\R}{\mathbb{R}}
\newcommand{\C}{\mathbb{C}}
\renewcommand{\H}{\mathcal{H}}
\newcommand{\B}{\mathcal{B}}
\newcommand{\A}{\mathcal{A}}
\newcommand{\K}{\mathcal{K}}
\newcommand{\M}{\mathcal{M}}
\newcommand{\vphi}{\varphi}
\newcommand{\Hess}{\operatorname{Hess}}

\newcommand{\J}{\mathscr{J}}
\newcommand{\D}{\mathscr{D}}
\renewcommand{\l}{\ell}
\newcommand{\Tr}{\text{Tr}}
\newcommand{\Var}{\text{Var}}
\newcommand{\HS}{\text{HS}}

\renewcommand{\P}{\mathcal{P}}

\newcommand{\ul}[1]{\underline{#1}}

\newcommand{\I}{\text{I}}
\newcommand{\II}{\text{II}}
\newcommand{\III}{\text{III}}

\newcommand{\<}{\left\langle}
\renewcommand{\>}{\right\rangle}
\renewcommand{\Re}[1]{\text{Re}\ #1}
\renewcommand{\Im}[1]{\text{Im}\ #1}
\newcommand{\dom}[1]{\text{dom}\,#1}
\renewcommand{\i}{\text{i}}
\renewcommand{\mod}[1]{(\operatorname{mod}#1)}
\newcommand{\mb}[1]{\mathbb{#1}}
\newcommand{\mc}[1]{\mathcal{#1}}
\newcommand{\mf}[1]{\mathfrak{#1}}
\newcommand{\im}{\operatorname{im}}

\newcommand{\TODO}[1]{{\color{red}\textbf{TODO: }{#1}}}


\newtheorem{thm}{Theorem}[section]
\newtheorem{prop}[thm]{Proposition}
\newtheorem{lem}[thm]{Lemma}
\newtheorem{cor}[thm]{Corollary}
\newtheorem{innercthm}{Theorem}
\newenvironment{cthm}[1]
  {\renewcommand\theinnercthm{#1}\innercthm}
  {\endinnercthm}
\newtheorem{innerclem}{Lemma}
\newenvironment{clem}[1]
  {\renewcommand\theinnerclem{#1}\innerclem}
  {\endinnerclem}

\theoremstyle{definition}
\newtheorem{defi}[thm]{Definition}
\newtheorem{ex}[thm]{Example}
\newtheorem{innercex}{Example}
\newenvironment{cex}[1]
	{\renewcommand\theinnercex{#1}\innercex}
	{\endinnercex}
\newtheorem*{exs}{Examples}
\newtheorem{rem}[thm]{Remark}
\newtheorem{innercdefi}{Definition}
\newenvironment{cdefi}[1]
  {\renewcommand\theinnercdefi{#1}\innercdefi}
  {\endinnercdefi}

\newtheorem*{imp}{}


\title[Free Stein kernels and an improvement of the free LSI]{Free Stein kernels and an improvement of the free logarithmic Sobolev inequality}

\author{Max Fathi$^\bullet$}
\address{UC Berkeley, Department of Mathematics}
\email{maxf@berkeley.edu}
\thanks{$\bullet$ Research supported by NSF FRG grant DMS-1361122}

\author{Brent Nelson$^\circ$}
\address{UC Berkeley, Department of Mathematics}
\email{brent@math.berkeley.edu}

\thanks{$\circ$ Research supported by NSF grants DMS-1502822}

\begin{abstract}

We introduce a free version of the Stein kernel, relative to a semicircular law. We use it to obtain a free counterpart of the HSI inequality of Ledoux, Peccatti and Nourdin, which is an improvement of the free logarithmic Sobolev inequality of Biane and Speicher, as well as a rate of convergence in the (multivariate) entropic free Central Limit Theorem. We also compute the free Stein kernels for several relevant families of self-adjoint operators. 
\end{abstract}

\maketitle


\section*{Introduction}

\subsection*{Entropy}

The classical notion of entropy was first introduced by Boltzmann in his work on the kinetic theory of gases, and since then has played an ubiquitous role in many fields, such as statistical physics, information theory and mathematics. The classical entropy takes the form
$$H(\nu) := \int{f \log f dx} \hspace{5mm} \nu = f \hspace{1mm} dx.$$
Following the usual convention in probability, this is the opposite of the physical entropy. Often, we are interested in a relative entropy with respect to a non-negative measure $\mu = e^{-V}dx$, which takes the form
$$H_{\mu}(\nu) = \int{\frac{d\nu}{d\mu}\log \left(\frac{d\nu}{d\mu}\right)d\mu} = H\left(\frac{d\nu}{dx}\right) + \int{Vd\nu}.$$
We would like to point out that physicists would call this quantity a free energy (the second term plays the role of an internal energy, and we omitted the temperature, since it will play no role in this work and can be absorbed in the potential). When the reference measure is the Lebesgue measure, the relative entropy coincides with the (simple) entropy. The relative entropy notably appears in Sanov's theorem as the large deviations rate function for sums of i.i.d. random variables. 

In \cite{VoiII}, Voiculescu introduced a notion of free entropy $\chi$ which plays the role of entropy in free probability theory. Most notably, it is monotone along the free central limit theorem, and the associated free energy/relative entropy appears as the large deviations rate function for the empirical spectrum of large random matrices. It also detects freeness: an $n$-tuple $(x_1,..,x_n)$ of non-commutative random variables with $\sum \chi(x_i)>-\infty$ is free iff it satisfies $\chi(x_1,..,x_n) = \sum \chi(x_i)$. 

The  microstates free entropy relative to a potential $V$ can then be defined as 
$$\chi(x_1,..,x_n | V) = \tau(V(x_1,..,x_n)) - \chi(x_1,..,x_n).$$
As in the classical case, we should think of this quantity as playing the role of a free energy. 

Several of the uses of the entropy in classical probability are consequences of its property of being monotone along the heat flow. However, for technical reasons, it is hard to study the analogous property for the microstates free entropy (which would be monotonicity along convolutions with a free semicircular random variable). It is not even clear whether the microstates free entropy can be differentiated along the free heat flow. To get around this difficulty, Voiculescu introduced a second notion of entropy, called the non-microstates free entropy. This second notion of entropy was specifically constructed so that the Fisher information could be easily introduced as the derivative of the entropy along the free heat flow. Since this work will mainly deal with logarithmic Sobolev inequalities, for which this property is essential, we shall mostly deal with this second notion of entropy. As is standard, the non-microstates free entropy (relative to $V$) shall be denoted by $\chi^*$ (resp. $\chi^*(x_1,\ldots, x_n|V)$). 

In dimension one, the two notions of entropy coincide, and have an explicit expression (as we shall see in Section \ref{section_1d}), but for $n$-tuples of non-commutative random variables with $n \geq 2$ this is an open problem.

\subsection*{Logarithmic Sobolev inequalities}

The classical logarithmic Sobolev inequality (LSI for short)
$$H_{\mu}(\nu) \leq \frac{1}{2\rho}\int{\frac{|\nabla f|^2}{f}d\mu}; \hspace{3mm} \nu = f\mu.$$
was introduced by Gross in \cite{Gr75}, where it was shown to be equivalent to hypercontractivity of the semigroup of a diffusion process naturally associated with $\mu$, a remarkable smoothing property that states that the action of the semigroup sends the space $L^p(\mu)$ into $L^q(\mu)$ for any $q > p > 1$ after a finite time. This inequality can be thought of as a kind of Sobolev inequality, where instead of embedding the Sobolev space $H^1$ into $L^p$ for some $p > 2$, we embed it into the Orlicz space $L^2\log L$. It has played an important role in many aspects of stochastic analysis: it is used to study long-time behavior for diffusion processes, it implies Gaussian concentration for the measure $\mu$ and has also been used to study large-scale behavior of interacting particle systems and large random matrices. It also implies Talagrand's transport-entropy inequality \cite{OV00}. We refer to \cite{ABC+} for more about its uses. 

A fundamental property of this inequality is that, unlike for the usual Sobolev inequalities, the dimension plays no role: when it holds for a metric-measure space, then it holds for products of the space, with the same constant. This property is the source of many interesting applications in infinite-dimensional analysis, such as dimension-free concentration bounds. 

When the reference measure is the standard Gaussian measure on $\R^n$, it is known that the LSI holds with sharp constant $\rho = 1$, and moreover cases of equality are explicitly known, and are standard gaussians with arbitrary mean \cites{Car91a,Car91b}. A classical theorem of Bakry and \'Emery \cite{BE85} states that if the reference measure on $\R^n$ is of the form $\mu = e^{-V}dx$ with $\Hess V \geq \rho Id$ for some $\rho > 0$, then it satisfies a LSI with constant $\rho$. 

In the free case, the LSI (if it holds) takes the form
$$E(X|V) - E(Y|V) \leq \frac{1}{2\rho}\Phi^*(X|V)$$
where $E\in \{\chi,\chi^*\}$, $V$ is some reference potential, $X$ is an $n$-tuple of free random variables, $Y$ is another $n$-tuple that minimizes the free entropy relative to $V$, and $\Phi^*(X|V)$ is the non-microstates free Fisher information of $X$ relative to the potential $V$. This inequality was obtained by Biane and Speicher \cite{BS01} when the reference potential is quadratic and $Y$ is therefore an $n$-tuple of free semicircular random variables.  

We shall still manage to obtain a LSI involving the microstates free entropy in the situation where the potential is uniformly convex. However, the improved versions we shall discuss in the next subsection rely on semigroup arguments, and we shall only establish them for the non-microstates free entropy. 

Unlike in the classical setting, we do not know whether the free LSI automatically implies the free Talagrand inequality (which has been studied in \cite{HPU04} for example), except in dimension one \cite{Led05}. The main obstacle is that the existence of free monotone transport is not yet well-understood and examples are quite rare.

\subsection*{Improved functional inequalities}

In the past decade, there has been a lot of interest in improved or quantitative functional inequalities. Informally, the situation is as follows: given a sharp functional inequality of the form $F(f) \leq G(f)$ for which the cases of equality are explicitly known, can one \emph{improve} the functional inequality by taking into account how far $f$ is from the set $\mathcal{F}$ of functions for which equality holds. Typical improvements that have been considered take the form
$$\delta(f) = G(f) - F(f) \geq c\ d(f, \mathcal{F})^{\alpha}$$
where $d$ is a well-chosen distance on the set of functions considered. Examples of inequalities for which such improvements have been obtained include isoperimetric inequalities \cite{FMP10}, the Brunn-Minkowski inequality \cite{EK14}, Sobolev inequalities \cite{CFMP09} and many others. In the last few years, there has been several partial results for the Gaussian logarithmic Sobolev inequality \cites{BGG15, BGRS14, FIL15, IM14, LNP16} but they are not entirely satisfactory. The problem of obtaining a deficit estimate in terms of a simple distance that is valid for all functions and behaves well in high dimension is still open at the time of writing. 

One of the results presented in this work is an adaptation of a result of \cite{BGRS14} on deficit estimates for the LSI to the free case. We show that for a $n$-tuple of non-commutative random variables with variance smaller than $n$, the deficit in the free LSI (with respect to a $n$-tuple of free semicircular laws) is bounded from below by a function of the free Fisher information. 

Another type of improvement to functional inequalities one can seek is to take into account another functional $H$ and seek a functional inequality of the form $F(f) \leq c(G(f), H(f))$ for some function $c$ satisfying $c(x, y) \leq x \hspace{1mm} \forall y$. Informally, this type of improvement consists of using extra information (the size of $H$) to learn more about how $F$ and $G$ relate. An improvement of this form for the Gaussian logarithmic Sobolev inequality, the so-called HSI inequality, has been obtained by Ledoux, Nourdin and Peccatti in \cite{LNP15}, involving the Stein deficit (see \ref{def_classical_stein}): 
	\begin{equation}
		H_{\gamma}(\nu) \leq \frac{S(\nu|\gamma)^2}{2}\log\left(1 + \frac{I_{\gamma}(\nu)}{S(\nu|\gamma)^2}\right)
	\end{equation}
for all centered measures $\nu$, where $H_{\gamma}$ stands for the relative entropy with respect to the Gaussian measure, and $I_{\gamma}$ for the Fisher information with respect to $\gamma$. Since $\log(1 + x) < x$ whenever $x > 0$, this is indeed a strict improvement of the classical LSI. This inequality has some applications to the study of long-time behavior of the Ornstein-Uhlenbeck semigroup, and can be used to yield rates of convergence in the entropic central limit theorem. 

One of the results of this work is to obtain a free analogue of this HSI inequality, involving a non-commutative analogue of the Stein discrepancy. We shall use this inequality to establish a rate of convergence in entropy for the multivariate free Central Limit Theorem. 

Other free functional inequalities for which improvements could be investigated include the one-dimensional Brunn-Minkowski inequality of \cite{Led05} and the one-dimensional transport-entropy inequality of \cite{HPU04}.

\subsection*{Acknowledgements}
We would like to thank Qiang Zeng for some helpful discussions we had in the early stages of this paper. We would also like to thank Dan-Virgil Voiculescu for his suggestions and encouragement.


\section{Definitions and Notation}

In this section, $M$ is a von Neumann algebra with a faithful normal state $\vphi$, and $X=(x_1,\ldots, x_n)\in M^n$ is an $n$-tuple of self-adjoint operators. We recall several free probabilistic quantities associated to $X$ and also define the quantities which appear in the improved free LSI.


\subsection{Norms}

Recall that via the GNS construction, $\vphi$ defines an inner product on $M$:
	\[
		\<x,y\>_\vphi:=\vphi(y^*x),\qquad x,y\in M.
	\]
The completion of $M$ with respect to the induced norm $\|\cdot\|_\vphi$ is denoted $L^2(M,\vphi)$. When the state $\vphi$ is clear from the context, we will abuse the notation and write $\|\cdot\|_2$ and $L^2(M)$.

For the von Neumann tensor product $M\bar\otimes M^{op}$ with faithful normal state $\vphi\otimes\vphi^{op}$, we can consider the Hilbert space $L^2(M\bar\otimes M^{op},\vphi\otimes\vphi^{op})$. This space can be identified with $\text{HS}(L^2(M))$ (the Hilbert--Schmidt operators on $L^2(M)$) via the following map
	\[
		x\otimes y^{op}\mapsto \<y,\cdot\,\>_2 x,\qquad x,y\in M.
	\]
In light of this identification, we will write $\<\cdot\,,\cdot\>_{\HS}$ for $\<\cdot\,,\cdot\>_{\vphi\otimes\vphi^{op}}$, when $\vphi$ is clear from the context.

For the $n$-tuple $X$ we define $\|X\|:=\max_j \|x_j\|$. We will also write $C^*(X)$ and $W^*(X)$ for the $C^*$-algebra and von Neumann algebra generated by $x_1,\ldots, x_n$, respectively. When $\vphi$ is clear, we write $\Var(X)$ for $\sum_{j=1}^n \vphi(x_j^*x_j)$. Given another $n$-tuple $Y=(y_1,\ldots, y_n)\in M^n$, we write
	\[
		\<X,Y\>_{\vphi}=\sum_{j=1}^n \<x_j, y_j\>_\vphi.
	\]
When $\vphi$ is clear we write $\<X,Y\>_2$. 

For $A, B\in L^2(M_n(M\bar\otimes M^{op}), (\vphi\otimes \vphi^{op})\circ \Tr)$ (the $n\times n$ matrices with entries in $L^2(M\bar\otimes M^{op},\vphi\otimes\vphi^{op})$), we write
	\[
		\<A,B\>_{(\vphi\otimes \vphi^{op})\circ \Tr}=\sum_{j,k=1}^n \<[A]_{jk}, [B]_{jk}\>_{\vphi\otimes \vphi^{op}}.
	\]
Again, when $\vphi$ is clear we simply write $\<A,B\>_{HS}$.


\subsection{Formal power series and non-commutative derivations}

We let $t_1,\ldots, t_n$ be non-commutating self-adjoint indeterminates, collected as the $n$-tuple $T=(t_1,\ldots, t_n)$. We let $\mathscr{P}$ denote the set of polynomials in these indeterminates. For a polynomial $p\in \mathscr{P}$ and a monomial $m$, we let $c_m(p)\in \C$ denote the coefficient of $m$ in $p$. After \cite{GS14}, for each $R>0$ we define the following norm
	\[
		\|p\|_R = \sum_m |c_m(p)| R^{\deg(m)},
	\]
where the (finite) sum is over all monomials appearing in $p$. The completion of $\mathscr{P}$ with respect to the $\|\cdot\|_R$-norm, a Banach algebra, is denoted $\mathscr{P}^{(R)}$ and can be thought of power series with radius of convergence at least $R$. Given an $n$-tuple of power series $F=(f_1,\ldots, f_n)\in (\mathscr{P}^{(R)})^n$, we define $\|F\|_R:=\max_j \|f_j\|_R$.

After \cite{VoiV}, for each $j=1,\ldots, n$ the \emph{$j$-th cyclic derivative} $\D_j$ is defined on monomials $m\in \mathscr{P}$ by
	\begin{align*}
		\D_j m = \sum_{m=at_jb} ba,
	\end{align*}
and extended linearly to all of $\mathscr{P}$. Then for $p\in \mathscr{P}$,
	\[
		\D p=(\D_1 p,\ldots, \D_n p)
	\]
is the \emph{cyclic derivative of $p$}. Similarly, for each $j=1,\ldots, n$, the \emph{$j$-th free difference quotient} $\partial$ is defined on monomials $m\in \mathscr{P}$ by
	\[
		\partial_j m = \sum_{m=at_j b} a\otimes b^{op},
	\]
and extended linearly to all of $\mathscr{P}$. For $P=(p_1,\ldots, p_n)\in \mathscr{P}^n$, we define
	\[
		\J P = \left(\begin{array}{cccc} \partial_1 p_1 & \partial_2 p_1 & \cdots & \partial_n p_1\\
									\vdots & \vdots & \ddots & \vdots \\
									\partial_1 p_n & \partial_n p_n & \cdots &\partial_np_n \end{array}\right) \in M_n(\mathscr{P}\otimes \mathscr{P}^{op}).
	\]
By \cite[Corollary 3.8]{NZ15}, $\J\D$ extends to $\mathscr{P}^{(R)}$, for any $R>0$, and is valued in the Banach algebra $M_n(\mathscr{P}\hat\otimes_R \mathscr{P}^{op})$. Here $\mathscr{P}\hat\otimes_R \mathscr{P}^{op}$ is the completion of the algebraic tensor product $\mathscr{P}^{(R)}\otimes (\mathscr{P}^{(R)})^{op}$ with respect to the projective tensor norm (see \cite[Section 3.1]{NZ15} for further details).

Given $p\in \mathscr{P}$, we will write $p(X)$ for the image of $p$ under the map defined by $t_{i_1}\cdots t_{i_d}\mapsto x_{i_1}\cdots x_{i_d}$. For any $R\geq \|X\|$, this map extends to $\mathscr{P}^{(R)}$, on which we use the same notation. We let $\mathscr{P}(X)^{(R)}$ denote the image of this extension, which we note lies in $C^*(X)$.

Similarly, given another $n$-tuple $Y=(y_1,\ldots,y_n)\in M^n$ and $\eta\in \mathscr{P}\otimes\mathscr{P}^{op}$, we write $\eta[X,Y]$ for the image of $\eta$ under the map $t_{i_1}\cdots t_{i_d}\otimes t_{j_1}\cdots t_{j_e}\mapsto x_{i_1}\cdots x_{i_d}\otimes y_{j_1}\cdots y_{j_e}$. When $Y=X$, we simply write $\eta[X]$. For $R\geq \|X\|,\|Y\|$, this map extends to $\mathscr{P}\hat\otimes_R\mathscr{P}^{op}$, on which we use the same notation. We let $\mathscr{P}(X)\hat\otimes_R\mathscr{P}(Y)$ denote the image of this extension, which we note lies in $C^*(X)\otimes_{\text{min}} C^*(Y)^{op}$


\subsection{Relative free Fisher information}

Recall from \cite{VoiV} that the \emph{conjugate variables to $X$} are vectors $\xi_{x_1},\ldots, \xi_{x_n}\in L^2(W^*(X),\vphi)$ satisfying for each $j=1,\ldots, n$ and each $p\in \mathscr{P}$
	\begin{align}\label{conj_var_equation}
		\<\xi_{x_j}, p(X)\>_2=\<1\otimes 1^{op}, [\partial_j p](X)\>_{\HS}.
	\end{align}
The \emph{free Fisher information for $x_1,\ldots, x_n$} is defined in \cite{VoiV} as the quantity
	\[
		\Phi^*(x_1,\ldots,x_n)=\sum_{j=1}^n \|\xi_{x_j}\|_2^2,
	\]
when the conjugate variables $\xi_{x_1},\ldots, \xi_{x_n}$ to $X$ exist, and as $+\infty$ otherwise. We will also denote this quantity by $\Phi^*(X)$.

\begin{defi}
Given $\mu\in \R$ and $\sigma>0$, by a \emph{free $(\mu,\sigma^2)$-semicircular $n$-tuple} we shall mean an $n$-tuple $(s_1,\ldots, s_n)$ of semicircular operators $s_1,\ldots, s_n$ that are freely independent from each other and satisfy $\vphi(s_j)=\mu$ and $\vphi(s_j^2)=\sigma^2$ for each $j=1,\ldots, n$.
\end{defi}

The conjugate variables to a free $(0,1)$-semicircular $n$-tuple $S=(s_1,\ldots, s_n)$ are simply $s_1,\ldots, s_n$. Hence $\Phi^*(S)=\Var(S)=n$. More generally, for $\rho>0$, the conjugate variables to $\frac{1}{\sqrt{\rho}}S$ (a free $(0,\rho^{-1})$-semicircular $n$-tuple) are $\sqrt\rho s_1,\ldots, \sqrt\rho s_n$ so that $\Phi^*(\frac{1}{\sqrt{\rho}} S)=n\rho$.

\begin{rem}\label{alg_free}
Recall that if the conjugate variables to $X$ exist and $R>\|X\|$, then by \cite[Theorem 2.5]{MSW15} and \cite[Lemma 37]{Dab14} the map $\mathscr{P}^{(R)}\ni f\mapsto f(X)$ is injective. Consequently, $\D_j$, $\D$, and $\partial_j$ induce derivations on $\mathscr{P}(X)^{(R)}$, which we denote $\D_{x_j}$, $\D_X$, and $\partial_{x_j}$, respectively. Similarly, $\J$ induces a derivation valued in $M_n(\mathscr{P}(X)\hat\otimes_R\mathscr{P}(X))$, which we denote by $\J_X$.
\end{rem}

For $V\in \mathscr{P}^{(R)}$, $R>\|X\|$, we say that the joint law of $X$ with respect to $\vphi$ is a \emph{free Gibbs state with potential $V$} if for each $j=1,\ldots, n$ and each $p\in \mathscr{P}$
	\[
		\<[\D_j V](X), p(X)\>_2=\<1\otimes 1^{op}, [\partial_j p](X)\>_{\HS}.
	\]
That is, if the conjugate variables to $X$ are given by $[\D_1 V](X),\ldots, [\D_n V](X)$. Equivalently, the following equation holds for all $P\in \mathscr{P}^n$:
	\[
		\<[\D V](X), P(X)\>_2 = \< (1\otimes 1^{op})\otimes I_n, [\J P](X)\>_{HS},
	\]
where $I_n\in M_n(\C)$ is the $n\times n$ identity matrix. For some simple examples, consider for each $\rho>0$ the potential
	\[
		V_\rho:=\frac{\rho}{2}\sum_{j=1}^n t_j^2 \in\mathscr{P}.
	\]
Then the joint law of a free $(0,\rho^{-1})$-semicircular $n$-tuple is a free Gibbs state with potential $V_\rho$.

\begin{defi}
For $X=(x_1,\ldots, x_n)\in M^n$, let $V\in \mathscr{P}^{(R)}$ for $R> \|X\|$. We define the \emph{free Fisher information of $x_1,\ldots, x_n$ relative to $V$} to be the quantity
	\[
		\Phi^*(x_1,\ldots, x_n\mid V):= \sum_{j=1}^n \|\xi_{x_j} - [\D_j V](X)\|_2^2
	\]
if the conjugate variables $\xi_{x_1},\ldots,\xi_{x_n}$ to $X$ exist, and as $+\infty$ otherwise. We also denote this quantity by $\Phi^*(X\mid V)$. In this context, we refer to $V$ as a \emph{potential}.
\end{defi}

The quantity $\Phi^*(X\mid V)$ is meant to measure how close the conjugate variables $\xi_{x_j}$ are to $[\D_j V](X)$, and so measures in some sense how close the joint law of $X$ is to a free Gibbs state with potential $V$.


\subsection{Relative non-microstates free entropy}

Let $S=(s_1,\ldots, s_n)$ be a free $(0,1)$-semicircular $n$-tuple, free from $x_1,\ldots, x_n$. Then the \emph{non-microstates free entropy of $x_1,\ldots, x_n$} is defined in \cite{VoiV} to be the quantity
	\[
		\chi^*(x_1,\ldots, x_n)=\frac{1}{2}\int_0^\infty \left(\frac{n}{1+t} - \Phi^*(X+\sqrt{t} S)\right)\ dt + \frac{n}{2}\log(2\pi e),
	\]
which we also denote by $\chi^*(X)$.

If $X$ is a free $(0,\rho^{-1})$-semicircular family, then $X+\sqrt{t} S$ is a free $(0,\rho^{-1}+t)$-semicircular family and therefore $\Phi^*(X+\sqrt{t} S)=n(\rho^{-1}+t)^{-1}$. From this it is easy to compute $\chi^*(X)=\frac{n}{2}\log(2\pi e\rho^{-1})$.

\begin{defi}
For $X=(x_1,\ldots, x_n)\in M^n$, let $V\in \mathscr{P}^{(R)}$ for $R\geq \|X\|$, we define the \emph{non-microstates free entropy of $x_1,\ldots, x_n$ relative to $V$} to be the quantity
	\[
		\chi^*(x_1,\ldots, x_n\mid V):=\vphi(V(X)) - \chi^*(x_1,\ldots, x_n),
	\]
which we will also denote by $\chi^*(X\mid V)$. In this context, we refer to $V$ as a \emph{potential}.
\end{defi}

Since $\chi^*(\,\cdot\,)$ is maximized (for fixed variance) by an $n$-tuple of free semicircular operators, it is easy to see that $\chi^*(\,\cdot\ \mid V_\rho)$ is minimized by a free $(0,\rho^{-1})$-semicircular $n$-tuple.


\subsection{Relative microstates free entropy}

For each $k\in \N$, let $M_{k}^{sa}$ denote the space of $k\times k$ self-adjoint matrices. For $l\in \N$ and $\epsilon>0$ consider the set
	\[
		\Gamma(X;k,l,\epsilon)=\{ Y\in (M_{k}^{sa})^n\colon |\frac{1}{n}\Tr(p(Y)) - 	\vphi(p(X))|<\epsilon,\ \forall p\in \mathscr{P} \text{ with } \text{deg}(p)\leq l\}.
	\]
The \emph{microstates free entropy of $x_1,\ldots, x_n$} is then defined in \cite{VoiII} as the quantity
	\[
		\chi(x_1,\ldots, x_n)=\inf_{\epsilon, l} \limsup_{k\to\infty} \frac{1}{k^2}\log(\text{Vol}(\Gamma(X;k,l,\epsilon))) + \frac{n}{2}\log{k},
	\]
which we will also denote by $\chi(X)$.

\begin{defi}
For $X=(x_1,\ldots, x_n)\in M^n$, let $V\in \mathscr{P}^{(R)}$ for $R\geq \|X\|$, we define the \emph{microstates free entropy of $x_1,\ldots, x_n$ relative to $V$} to be the quantity
	\[
		\chi(x_1,\ldots, x_n\mid V):=\vphi(V(X)) - \chi(x_1,\ldots, x_n),
	\]
which we will also denote by $\chi(X\mid V)$. In this context, we refer to $V$ as a \emph{potential}.
\end{defi}

The functional $\chi(\,\cdot\mid V)$ has been previously considered in \cite[Section 3.7]{Voi02} and \cite[Theorem 5.1]{GS14}, where the following lemma was observed.

\begin{lem}\label{minimizer_is_free_Gibbs_state}
Let $\vphi=\tau$ be a trace. For $V=V^*\in \mathscr{P}^{(R)}$, suppose $Y=(y_1,\ldots, y_n)\in M^n$ satisfies $\|Y\|<R$ and minimizes $\chi(\,\cdot\mid V)$. Then joint law of $Y$ is a free Gibbs state with potential $V$.
\end{lem}
\begin{proof}
Let $P=(p_1,\ldots, p_n)\in \mathscr{P}^n$ be an $n$-tuple of polynomials. By \cite[Proposition 1.3]{VoiIV}, we have
	\begin{align*}
		\frac{d}{dt} \chi(y_1+tp_1(Y),\ldots, y_n+tp_n(Y))\mid_{t=0} = \sum_{j=1}^n \tau\otimes\tau^{op}([\partial_j p_j](Y)) = (\tau\otimes\tau^{op})\circ\Tr( [\J P](Y)).
	\end{align*}
It is also easy to compute that
	\begin{align*}
		\frac{d}{dt} \tau( V(y_1+tp_1(Y),\ldots, y_n+tp_n(Y)))\mid_{t=0} = \sum_{j=1}^n \tau( [\partial_j V](Y)\# p_j(Y))=\tau( [\D V](Y)\cdot P(Y)),
	\end{align*}
where $(a\otimes b^{op})\# c=acb$. So, by virtue of minimizing the microstates free entropy relative to $V$ we have
	\begin{align*}
		0 = \frac{d}{dt}\chi(y_1+tp_1(Y),\ldots, y_n+tp_n(Y)\mid V)\mid_{t=0} =  \tau( [\D V](Y)\cdot P(Y)) - (\tau\otimes\tau^{op})\circ\Tr( [\J P](Y));
	\end{align*}
that is, the joint law of $y_1,\ldots, y_n$ is a free Gibbs state with potential $V$.
\end{proof}

%


\subsection{Free Stein Discrepancy}

In this section, we introduce a free analogue of the Stein discrepancy. In commutative probability, the Stein discrepancy is a way to measure how far away a probability measure is from a given reference measure by looking at how badly it violates a set of relations that characterize the reference measure. These relations are often defined via an integration by parts formula the reference measure satisfies. For the Gaussian measure, the formula usually employed is
\begin{equation} \label{eq_stein_gamma}
\int{x \cdot \nabla f d\gamma} = \int{\Delta f d\gamma} \hspace{5mm} \forall f \in C^{\infty}_c(\R^n).
\end{equation}
In general, one finds such an integration by parts formula by characterizing the reference measure $\mu$ as the unique invariant probability measure of a well-chosen reversible Markov process with generator $\mathcal{L}$, and then using the relation $\int{\mathcal{L}f d\mu} = 0 \hspace{5mm} \forall f$. For the Gaussian measure, we get the above relation when using the Ornstein-Uhlenbeck flow as such a Markov process. 

\begin{defi} \label{def_classical_stein}
Let $\nu$ be a centered probability measure on $\R^n$. The Stein kernel of $\nu$ with respect to the Gaussian measure $\gamma$ on $\R^n$ is the matrix-valued function $\tau : \R^n \longrightarrow S_n(\R)$ such that
	\begin{align}\label{comm_Stein_kernel_eq}
		\int{x\cdot \nabla f d\nu} = \int{\langle \tau, \Hess f \rangle_{HS} d\nu}.
	\end{align}
The Stein discrepancy of $\nu$ with respect to $\gamma$ is then given by
$$S(\nu|\gamma) := \left(\int{||\tau - Id||_{HS}^2d\nu}\right)^{1/2}$$
if the Stein kernel exists, and $+\infty$ if not. 
\end{defi}

The Stein kernel is a way to reformulate the integration by parts formula. It is easy to see that a probability measure satisfies \eqref{eq_stein_gamma} for all $f$ if and only if its Stein kernel is the identity. We also note that for a Stein kernel with respect to the Gaussian to exist, $\nu$ must be centered. We refer to \cite{LNP15} and the references therein for more about the Stein method in classical probability. 

We now introduce the non-commutative analogues of the Stein kernel and Stein discrepancy: 

\begin{defi}
For $X=(x_1,\ldots, x_n)\in M^n$, let $V\in\mathscr{P}^{(R)}$ for $R> \|X\|$. An element
	\[
		A\in L^2(M_n(M\bar\otimes M^{op}),(\vphi\otimes\vphi^{op})\circ\Tr)
	\]
is said to be a \emph{free Stein kernel for $X$ relative to $V$ (with respect to $\varphi$)} if
	\begin{align}\label{Stein_kernel_eqn}
		\<[\D V](X), P(X)\>_\vphi=\<A, [\J P](X)\>_{(\vphi\otimes\vphi^{op})\circ\Tr},\qquad \forall P\in \mathscr{P}^n.
	\end{align}
The \emph{free Stein discrepancy of $X$ relative to $V$} is defined as the quantity
	\[
		\Sigma^*(X\mid V)=\inf_A \|A- (1\otimes 1^{op})\otimes I_n\|_{(\varphi\otimes\varphi^{op})\circ\Tr},
	\]
where the infimum is taken over the set of Stein kernels of $X$ relative to $V$.
\end{defi}

Observe that (\ref{Stein_kernel_eqn}) is not the exact analogue of (\ref{comm_Stein_kernel_eq}). Indeed, the exact analogue would be
	\[
		\<[\D V](X), [\D p](X)\>_\vphi = \< A, [\J\D p](X)\>_{(\vphi\otimes\vphi^{op})\circ\Tr},\qquad \forall p\in \mathscr{P},
	\]
which is a weaker condition. Fortunately, all of the examples in Section \ref{examples} satisfy our stronger condition and therefore offer applications of the free HSI. We use this stronger condition in the proof of Lemma \ref{Stein_decay}, but the weaker condition would suffice if a conjecture of Voiculescu receives even a partially affirmative answer (see Remark \ref{not_exact_analogue}).

We note that since $[\J(\mathscr{P}^n)](X)$ need not be dense in $L^2(M_n(M\bar\otimes M^{op}),(\vphi\otimes\vphi^{op})\circ\Tr)$, the free Stein kernel may not be unique. Nevertheless, so long as at least one free Stein kernel exists, then its projection onto the $L^2$-closure of $[\J(\mathscr{P}^n)](X)$ is also a free Stein kernel which attains the infimum in the free Stein discrepancy. Furthermore, even the Stein kernel of a probability measure $\nu$ on $\R^d$ need not be unique for $d>1$: the space of Hessians is not dense.

When $\Phi^*(X)<+\infty$, Remark \ref{alg_free} implies that (\ref{Stein_kernel_eqn}) is equivalent to $A\in \dom(\J_X^*)$ with $\J_X^*(A)=[\mathscr{D}V](X)$. Also note that when the potential $V$ is self-adjoint and $\varphi$ is a trace, equation (\ref{Stein_kernel_eqn}) implies $A^\dagger=A$, where $[A^\dagger]_{ij}=[A]_{ij}^\dagger$ and $(a\otimes b^{op})^\dagger=b^*\otimes (a^*)^{op}$ for $a,b\in M$.

Just as with the free Fisher information relative to some potential $V$, the free Stein discrepancy relative to $V$ measures in some sense how close the $n$-tuple is to having as its joint law a free Gibbs state with potential $V$. Indeed, by the definition of a free Gibbs state, $X$ has such a joint distribution if and only if $(1\otimes 1^{op})\otimes I_n$ is a free Stein kernel for $X$ and consequently $\Sigma^*(X\mid V)=0$. In particular, we have $\Sigma^*(\frac{1}{\sqrt{\rho}} S\mid V_\rho)=0$. We will denote $(1\otimes 1^{op})\otimes I_n$ by $1$ when the context precludes any confusion.


\section{An Improved (Non-Microstates) Free log-Sobolev Inequality}

We now fix a von Neumann algebra $M$ with a faithful normal trace $\tau$ and an $n$-tuple $X=(x_1,\ldots, x_n)$ of self-adjoint operators. Provided they exist, the conjugate variables to $X$ will always be denoted $\xi_1,\ldots, \xi_n$. Let $S=(s_1,\ldots, s_n)\in M^n$ be a free $(0,1)$-semicircular $n$-tuple, free from $x_1,\ldots, x_n$. Fix $\rho>0$ and recall the potential
	\[
		V_\rho=\frac{\rho}{2}\sum_{j=1}^n t_j^2.
	\]
For each $t\geq 0$ and $j=1,\ldots, n$, define $x_j(t):=e^{-t} x_j + \sqrt{ 1- e^{-2t}} \frac{1}{\sqrt{\rho}} s_j$, collected as the $n$-tuple $X(t)=(x_1(t),\ldots, x_n(t))$. Let $\mc{E}_t\colon M\to W^*(X(t))$ denote the conditional expectation. By \cite[Proposition 3.7]{VoiV}, the conjugate variables to $X(t)$ always exist, and if we denote them by $\xi_1(t),\ldots, \xi_n(t)$ then for each $j=1,\ldots, n$
	\[
		\xi_j(t)=\frac{\sqrt{\rho}}{\sqrt{1-e^{-2t}}} \mc{E}_t(s_j).
	\]
Moreover, if the conjugate variables to $X$ exist then we have
	\[
		\xi_j(t)=e^t\mc{E}_t(\xi_j).
	\]
We collect the conjugate variables to $X(t)$ as the $n$-tuple $\Xi(t)=(\xi_1(t),\ldots, \xi_n(t))$.


\subsection{De Bruijn's formula and a free log-Sobolev inequality}

The following lemma gives the non-commutative version of de Bruijn's formula.

\begin{lem}\label{deBruijn}
For $\rho>0$ and $X(t)=(x_1(t),\ldots, x_n(t))$ as above,
	\[
		\frac{d}{dt} \chi^*(X(t)\mid V_\rho)= -\frac{1}{\rho} \Phi^*(X(t)\mid V_\rho).
	\]
\end{lem}
\begin{proof}
We first note that because $\D_j V_\rho = \rho t_j$, we have
	\begin{align*}
		\Phi^*(X(t)\mid V_\rho) &= \sum_{j=1}^n \|\xi_j(t) - \rho x_j(t)\|_2^2\\
			&= \sum_{j=1}^n \|\xi_j(t)\|_2^2 - 2\rho + \rho^2( e^{-2t}\|x_j\|_2^2 + (1-e^{-2t})\rho^{-1})\\
			&= \Phi^*(X(t)) - n\rho + \rho^2 e^{-2t}\Var(X) - n\rho e^{-2t}.
	\end{align*}
Next we compute
	\begin{align*}
		\frac{d}{dt} \tau(V_\rho(X(t))= \frac{d}{dt}\left[ \frac{\rho}{2} e^{-2t} \Var(X) + \frac{n}{2}(1-e^{-2t})\right] = -\rho e^{-2t}\Var(X) + n e^{-2t}.
	\end{align*}
By \cite[Proposition 7.8]{VoiV},
	\begin{align*}
		\chi^*(X(t)) = -nt + \chi^*\left( X+ \sqrt{\rho^{-1} (e^{2t} - 1)} S\right).
	\end{align*}
Let $r(t) = \rho^{-1}(e^{2t} - 1)$. Then by \cite[Proposition 7.5]{VoiV}  we have
	\begin{align*}
		\frac{d}{dt}\chi^*( X(t) ) &= -n + \frac{1}{2} \Phi^*(X+ \sqrt{r(t)}S) \frac{d r(t)}{dt}\\
			&= - n + \rho^{-1} e^{2t} \Phi^*\left(X+\sqrt{e^{2t} - 1}\frac{1}{\sqrt{\rho}}S\right)\\
			&= - n +  \rho^{-1} \Phi^*(X(t)).
	\end{align*}
Thus the claimed equality holds.
\end{proof}

\begin{lem}\label{exp_decay}
For $\rho>0$ and $X(t)=(x_1(t),\ldots, x_n(t))$ as above,
	\[
		\Phi^*(X(t)\mid V_\rho) \leq e^{-2t} \Phi^*(X\mid V_\rho).
	\]
\end{lem}
\begin{proof}
If $\Phi^*(X)=+\infty$, then there is nothing to prove. So assume $\Phi^*(X)<+\infty$ and let $\xi_1,\ldots, \xi_n$ denote the conjugate variables. Recall that for each $j=1,\ldots, n$
	\[
		\xi_j(t)  = e^t \mc{E}_t(\xi_j) = \frac{\sqrt{\rho}}{\sqrt{1-e^{-2t}}} \mc{E}_t (s_j).
	\]
As observed in Lemma \ref{deBruijn},
	\[
		\Phi^*(X(t)\mid V_\alpha) = \sum_{j=1}^n \|\xi_j(t)\|_2^2 -\rho (1+e^{-2t}) + \rho^2 e^{-2t}\|x_j\|_2^2
	\]
Now, for each $j=1,\ldots, n$ we have
	\begin{align*}
		\|\xi_j(t)\|_2^2 &= \left\| e^{-2t} e^t\mc{E}_t(\xi_j) + (1-e^{-2t}) \frac{\sqrt{\rho}}{\sqrt{1-e^{-2t}}} \mc{E}_t(s_j)\right\|_2^2\\
			&= \left\| \mc{E}_t\left( e^{-t} \xi_j + \sqrt{\rho(1-e^{-2t})} s_j\right)\right\|_2^2\\
			&\leq \left\| e^{-t} \xi_j + \sqrt{\rho(1-e^{-2t})} s_j\right\|_2^2\\
			&= e^{-2t} \|\xi_j\|^2 + \rho(1-e^{-2t}).
	\end{align*}
So, resuming our computation from above we see
	\begin{align*}
		\Phi^*(X(t)\mid V_\rho) &\leq \sum_{j=1}^n e^{-2t} \|\xi_j\|^2 + \rho (1-e^{-2t}) - \rho(1+e^{-2t}) + \rho^2 e^{-2t}\|x_j\|_2^2\\
			&= e^{-2t} \left(\sum_{j=1}^n \|\xi_j\|_2^2 - 2\rho + \rho^2\|x_j\|_2^2\right)\\
			&= e^{-2t} \Phi^*(X\mid V_\rho),
	\end{align*}
as claimed.
\end{proof}

From these two lemmas we obtain the non-commutative analogue of the logarithmic-Sobolev inequality from classical probability. In particular, this yields an alternate proof for an inequality from \cite[Section 7.2]{BS01} stated here in the multi-variable case.

\begin{cor}\label{log-Sob}
For $X$ an $n$-tuple of non-commutative random variables, $S$ a free $(0,1)$-semicircular $n$-tuple that is free from $X$, and any $\rho>0$
	\[
		\chi^*(X\mid V_\rho) - \chi^*(\sqrt{\rho}^{-1} S\mid V_\rho) \leq \frac{1}{2\rho} \Phi^*(X\mid V_\rho).
	\]
In particular, if $\Var(X)=\frac{n}{\rho}$ then
	\[
		\frac{n}{2} + \frac{n}{2} \log\left(\frac{2\pi e}{\rho} \right) - \chi^*(X) \leq \frac{1}{2\rho}\Phi^*(X).
	\]
\end{cor}
\begin{proof}
Using Lemmas \ref{deBruijn} and \ref{exp_decay} we have
	\begin{align*}
		\chi^*(X\mid V_\rho) - \chi^*(\sqrt{\rho}^{-1}S\mid V_\rho) &= \chi^*(X(0)\mid V_\rho) - \lim_{t\to\infty} \chi^*( X(t)\mid V_\rho) \\
			&= \int_0^\infty \rho^{-1}\Phi^*(X(t)\mid V_\rho)\ dt \\
			&\leq \rho^{-1}\Phi^*(X\mid V_\rho) \int_0^\infty e^{-2t}\ dt = \frac{1}{2\rho} \Phi^*(X\mid V_\rho).
	\end{align*}
The second inequality follows by expanding the definitions of the relative free entropy and Fisher information (and using $\chi^*(\sqrt{\rho}^{-1} S)= \frac{n}{2}\log\left(\frac{2\pi e}{\rho}\right)$).
\end{proof}


\subsection{Improving the free log-Sobelev Inequality via Stein's Method}

\begin{lem}\label{Stein_decay}
Suppose $\Phi^*(X\mid V_\rho)<+\infty$. For $\rho>0$ and $X(t)=(x_1(t),\ldots, x_n(t))$ as above
	\[
		\frac{1}{\rho}\Phi^*(X(t)\mid V_\rho) \leq \frac{e^{-4t}}{1-e^{-2t}} \Sigma^*(X\mid V_\rho)^2
	\]
\end{lem}
\begin{proof}
If $\Sigma^*(X\mid V_\rho)=+\infty$, there is nothing to prove. So assume $\Sigma^*(X\mid V_\rho)<+\infty$, and let $A\in L^2(M_n(W^*(X)\bar\otimes W^*(X)^{op}))$ be a free Stein kernel for $X$ relative to $V_\rho$ such that 
	\[
		\|A-1\|_{HS}=\Sigma^*(X\mid V_\rho).
	\]
Recall that since $\Phi^*(X)<+\infty$, this implies $\J_X^*(A)=\rho X$. Let $\xi_1,\ldots, \xi_n\in L^2(W^*(X))$ be the conjugate variables of $X$, collected in a vector as $\Xi=(\xi_1,\ldots, \xi_n)$ so that $\J_X^*(1)=\Xi$. Now, we compute:
	\begin{align}\label{strong_step}
		\Phi^*(X(t)\mid V_\rho) &= \sum_{j=1}^n \|\xi_j(t) - \rho x_j(t)\|_2^2\nonumber\\
			&= \sum_{j=1}^n \< e^{-t} \xi_j + \sqrt{\rho(1-e^{-2t})} s_j - \rho x_j(t) , \xi_j(t)-\rho x_j(t)\>_2\nonumber\\
			&=  e^{-t} \sum_{j=1}^n  \< \xi_j - \rho x_j, \xi_j(t) -  \rho x_j(t)\>_2\nonumber\\
			&= e^{-t} \< \Xi - \rho X,  \Xi(t) -\rho X(t)\>_2\nonumber\\
			&= e^{-t} \< \J_X^*[1 -  A],  \Xi(t)- \rho X(t)\>_2.
	\end{align}
Recall that on $\C\<x_1(t),\ldots, x_n(t)\>$
	\[
		e^t\J_X=\J_{X(t)}= \frac{\sqrt{\rho}}{\sqrt{1-e^{-2t}}} \J_S,
	\]
so that for any $Z\in \mc{E}_t(M)^n$ we have
	\[
		\<\J_X^*[1-A], Z\>_2 = \frac{e^{-t}\sqrt{\rho}}{\sqrt{1-e^{-2t}}} \< \J_S^*[1- A],Z\>_2.
	\]
So continuing our computation we have
	\begin{align*}
		\Phi^*(X(t)\mid V_\rho) &= \frac{e^{-2t}\sqrt{\rho}}{\sqrt{1-e^{-2t}}} \< \J_S^*[1 - A], \Xi(t) - \rho X(t)\>_2\\
			&\leq \frac{e^{-2t}\sqrt{\rho}}{\sqrt{1-e^{-2t}}} \|\J_S^*[1 - A]\|_2 \Phi^*(X(t)\mid V_\rho)^{1/2}.
	\end{align*}
Now, since $\partial_{s_j}\mid_{W^*(X)}\equiv 0$ for each $j=1,\ldots, n$ we have by \cite[Proposition 4.1]{VoiV}
	\begin{align*}
		\J_S^*[1-A] &= \left( \sum_{j=1}^n \partial_{s_j}^*[ \delta_{i=j} 1\otimes 1 - A_{ij}]\right)_{i=1}^n= \left(  s_i - \sum_{j=1}^n A_{ij}\# s_j \right)_{i=1}^n = [1-A]\# S.
	\end{align*}
Since $s_1,\ldots, s_n$ are free from $W^*(X)$, it is easy to see that
	\[
		\|[1-A]\#S\|_2 = \| 1 - A\|_{HS} = \Sigma^*(X\mid V_\rho)
	\]
Thus we have
	\[
		\Phi^*(X(t)\mid V_\rho )\leq \frac{e^{-2t}\sqrt{\rho}}{\sqrt{1-e^{-2t}}} \Sigma^*(X\mid V_\rho) \Phi^*(X(t)\mid V_\rho)^{1/2},
	\]
which is equivalent to the desired inequality.
\end{proof}

\begin{rem}\label{not_exact_analogue}
The equality on line (\ref{strong_step}) is precisely where we use the stronger version of the definition of a free Stein kernel. It is a conjecture of Voiculescu that conjugate variables always lie in the $L^2$-closure of cyclic gradients. If this even weakly true for $\Xi(t)$, then the equality on line (\ref{strong_step}) will hold for the weaker version of the definition of a free Stein kernel.
\end{rem}

With the previous lemmas in hand, the following proof is, aside from notation, identical to that of \cite[Theorem 2.2]{LNP15}. We still present it here, for the sake of convenience.

\begin{thm}
Let $X=(x_1,\ldots, x_n)$ be an $n$-tuple of non-commutative random variables, and $S=(s_1,\ldots, s_n)$ a free $(0,1)$-semicircular $n$-tuple, free from $X$. Then for any $\rho>0$ we have
	\[
		\chi^*(X\mid V_\rho) - \chi^*(\sqrt{\rho}^{-1}S\mid V_\rho) \leq \frac{1}{2} \Sigma^*(X\mid V_\rho)^2\log\left(1+\frac{\Phi^*(X\mid V_\rho)}{\rho \Sigma^*(X\mid V_\rho)^2}\right).
	\]
\end{thm}
\begin{proof}
If $\Phi^*(X\mid V_\rho)=+\infty$, then there is nothing to prove, so we assume that $\Phi^*(X\mid V_\rho)<+\infty$. If $\Sigma^*(X\mid V_\rho)=+\infty$, then the right-hand side reduces to $\frac{1}{2\rho}\Phi^*(X\mid V_\rho)$ and the inequality holds by Corollary \ref{log-Sob}. Thus we may assume $\Sigma^*(X\mid V_\rho)<+\infty$.

Using Lemmas \ref{deBruijn}, \ref{exp_decay}, and \ref{Stein_decay} we have
	\begin{align*}
		\chi^*(X\mid V_\rho) - \chi^*(\sqrt{\rho}^{-1}S\mid V_\rho) &=\int_0^\infty \rho^{-1} \Phi^*(X(t)\mid V_\rho)\ dt \\
			&\leq \rho^{-1} \Phi^*(X\mid V_\rho)\int_0^u e^{-2t}\ dt + \Sigma^*(X\mid V_\rho)^2 \int_u^{\infty} \frac{e^{-4t}}{1-e^{-2t}}\ dt\\
			&\leq \frac{1}{2\rho}\Phi^*(X\mid V_\rho)(1-e^{-2u}) + \frac{1}{2} \Sigma^*(X\mid V_\rho)^2(-e^{-2u} - \log(1-e^{-2u})).
	\end{align*}
Optimizing in $u$ yields $1-e^{-2u}=\frac{\rho\Sigma^*(X\mid V_\rho)^2}{\Phi^*(X\mid V_\rho) + \rho \Sigma^*(X\mid V_\rho)^2}$ and the desired inequality.
\end{proof}


\section{Examples of free Stein kernels}\label{examples}

In this section we present several examples $n$-tuples of self-adjoint operators with known free Stein kernels relative to the potential $V_1$, for which equation (\ref{Stein_kernel_eqn}) reduces to
	\[
		\<X, P\>_{\vphi}=\<A, \J_X P\>_{(\vphi\otimes\vphi^{op})\circ\Tr},\qquad \forall P\in \mathscr{P}(X)^n,
	\]
which, in turn, is equivalent to the system of equations
	\[
		\<x_j, p(X)\>_{\vphi} = \sum_{k=1}^n \<[A]_{jk}, [\partial_k p](X)\>_{\vphi\otimes\vphi^{op}},\qquad p\in\mathscr{P},
	\]
for $j=1,\ldots, n$.

Notably, these examples represent all known applications of the free monotone transport theorem of Guionnet and Shlyakhtenko \cite{GS14} (excepting the infinite variable case considered in \cite{NZ15}). That is, under certain hypothesis, the $n$-tuples of operators considered in these examples have been shown by free transport arguments to generate a von Neumann algebra isomorphic to the free group factor $L(\mathbb{F}_n)$ (the von Neumann algebra generated by a free $(0,1)$-semicircular $n$-tuple), or (in Example 3) a free Araki-Woods factor (the non-tracial analogue of $L(\mathbb{F}_n)$). Moreover, the free Stein discrepancy being small is a necessary condition under the hypotheses used in these free transport arguments. Indeed, the free Stein discrepancy is dominated by a norm derived from the projective tensor norm on $\mathscr{P}^{(R)}\otimes(\mathscr{P}^{(R)})^{op}$. At the core of the proofs of each of these applications is the fact that this latter norm must be sufficiently small (so that, for example, the free Stein kernel is invertible in a particular Banach algebra). In light of this, we posit the following conjecture. See \cite{GS14} for further details on free monotone transport.

\newtheorem*{conj}{Conjecture}
\begin{conj}
Suppose $X=(x_1,\ldots, x_n)$ has joint law $\mu_X$, and $A=A^\dagger$ is a free Stein kernel for $X$ relative to $V_1$. Then there exists $\epsilon>0$ such that if $\Sigma^*(X\mid V_1)<\epsilon$, then there exists free transport from the free semicircle law to $\mu_X$. In particular, we can embed $W^*(X)\hookrightarrow L(\mathbb{F}_n)$.
\end{conj}

Each of the following examples arises from a similar construction: Let $\mc{H}_\R$ let be a real Hilbert space with orthonormal basis $\{e_1,\ldots, e_n\}$, contained in some complex Hilbert space $\mc{H}$. Then the operators in our examples will be represented on some Hilbert space completion $\mc{F}$ of the following space:
	\[
		\mc{F}^{\text{finite}}:=\C\Omega\oplus\bigoplus_{d\geq 1} \mc{H}^{\otimes d},
	\]
where $\Omega$ is the vacuum vector (to be regarded as the generator of ``$\mc{H}^{\otimes 0}$''). We consider linear operators $\ell(e_1),\ldots, \ell(e_n)$ on $\mc{F}^{\text{finite}}$ defined by
	\begin{align*}
		\ell(e_j)\Omega &= e_j\\
		\ell(e_j) f_1\otimes\cdots\otimes f_d &= e_j\otimes f_1\otimes\cdots\otimes f_d,\qquad f_1,\ldots, f_d\in \mc{H}.
	\end{align*}
These operators, which will be bounded operators on all the completions of $\mc{F}^{\text{finite}}$, are known as \emph{left creation operators}, possibly with some additional adjective depending on the completion. Furthermore, the particular completion considered will determine the behavior of the adjoints $\ell(e_1)^*,\ldots, \ell(e_n)^*$, which are known as \emph{left annihilation operators} (possibly with some additional adjectives). In fact, for the sake of these examples this behavior is more enlightening than an explicit description of the inner product on $\mc{F}^{\text{finite}}$ which determines the completion $\mc{F}$, and so we forgo the latter.

In the following examples, the self-adjoint operators of interest are $x_j=\ell(e_j)+\ell(e_j)^*$, $j=1,\ldots, n$. We consider their distribution with respect to the vacuum state $\vphi(\,\cdot\,)=\<\,\cdot\ \Omega, \Omega\>_\mc{F}$. In particular, we will describe free Stein kernels for the $n$-tuple $X=(x_1,\ldots, x_n)$ relative to $V_1$.


\begin{cex}{1: $q$-semicircular systems \cite{BS91}}
Let $\mc{H}_\R=\R^n$, contained in $\mc{H}=\C^n$. For $q\in[-1,1]$, one can define an inner product on $\mc{F}^{\text{finite}}$ (see \cite[Lemma 3]{BS91} for an explicit formula) so that for each $j=1,\ldots, n$
	\begin{align*}
		\ell(e_j)^*\Omega &=0\\
		\ell(e_j)^* f_1\otimes\cdots\otimes f_d &= \sum_{k=1}^d q^{k-1}\<f_k,e_j\> f_1\otimes\cdots\otimes \hat{f}_k\otimes\cdots\otimes f_d,\qquad f_1,\ldots, f_d\in \C^n,
	\end{align*}
where $\hat{f}_k$ denotes omission from the tensor product. In this case, $\{\ell(e_j)\}_{j=1}^n$, $\{\ell(e_j)^*\}_{j=1}^n$, and $\{x_j\}_{j=1}^n$ are respectively known as \emph{left $q$-creation operators}, \emph{left $q$-annihilation operators}, and \emph{$q$-semicircular operators}. Also, the vacuum state is a trace which we denote $\tau$. When $q=0$, $X$ is a free $(0,1)$-semicircular $n$-tuple.

Towards exhibiting a free Stein kernel, for each $d\geq 1$ let $\pi_d\in \mc{B}(\mc{F})$ be the projection onto $(\C^n)^{\otimes d}$, and let $\pi_0\in\mc{B}(\mc{F}_q(n))$ be the projection onto $\C\Omega$. Consider the operator
	\[
		\Xi_q=\sum_{d\geq 0} q^d \pi_d.
	\]
When $q^2n<1$, $\Xi_q$ is a Hilbert--Schmidt operator on $\mc{F}$. Consequently, using the identification
	\begin{align*}
		L^2(W^*(X)\bar\otimes W^*(X)^{op}) &\cong \text{HS}(\mc{F})\\
			a\otimes b^\circ &\mapsto \<\,\cdot\ ,b^*\Omega\>_{\mc{F}} a\Omega,
	\end{align*}
we can consider $A=\Xi_q\otimes I_n\in L^2(M_n(W^*(X)\bar\otimes W^*(X)^{op}))$. It was observed by Shlyakhtenko in \cite{Shl04} that $\J_X^*(A)=X$ (see also \cite[Lemma 10]{Shl09}). Hence $A$ is a free Stein kernel for $X$ with respect to $V_1$, and
	\begin{align*}
		\Sigma^*(X\mid V_1)\leq \sqrt{n} \|\Xi_q - 1\otimes 1^\circ\|_{\tau\otimes \tau^{op}} = \frac{|q|n}{\sqrt{1-q^2 n}}.
	\end{align*}

\end{cex}


\begin{cex}{2: mixed $q$-semicircular systems \cite{BS94}}
As before, we let $\mc{H}_\R=\R^n$ and $\mc{H}=\C^n$. Given a symmetric matrix $Q=(q_{ij})\in M_n([-1,1])$, one can define an inner product on $\mc{F}^{\text{finite}}$ (see \cite[Section 3]{BS94} for a formula) so that for each $j=1,\ldots n$
	\begin{align*}
		\ell(e_j)^*\Omega &=0\\
		\ell(e_j)^* e_{i_1}\otimes \cdots \otimes e_{i_d} &= \sum_{k=1}^n q_{ji_1}\cdots q_{ji_{k-1}} \<e_{i_k}, e_j\> e_{i_1}\otimes \cdots\otimes \hat{e}_{i_k}\otimes \cdots \otimes e_{i_d}.
	\end{align*}
The operators $x_j=\ell(e_j)+\ell(e_j)^*$ are called \emph{mixed $q$-Gaussian variables}, and the vacuum state is again a trace, which we denote by $\tau$.

As before, a free Stein kernel will arise from Hilbert--Schmidt operators. However, in order to properly describe them we require a finer orthogonal decomposition of $\mc{F}$. For $d\geq 1$ and $1\leq i_1\leq i_2\leq\cdots \leq i_d\leq n$ define
	\[
		\mc{F}_{(i_1,\ldots, i_n)} = \text{span}\{ e_{i_{\sigma(1)}}\otimes\cdots\otimes e_{i_{\sigma(d)}}\colon \sigma\in S_d\},
	\]
where $S_d$ is the permutation group on $d$ elements. Then these subspaces along with $\C\Omega$ offer an orthogonal decomposition of $\mc{F}$. Let $\pi_{(i_1,\ldots, i_d)}$ be the orthogonal projection onto $\mc{F}_{(i_1,\ldots, i_d)}$ and let $\pi_0$ be the orthogonal projection onto $\C\Omega$. For each $j=1,\ldots, n$, define
	\[
		\Xi_j=\pi_0+\sum_{d\geq 1}\sum_{1\leq i_1\leq \cdots \leq i_d\leq n} q_{ji_1}\cdots q_{ji_d} \pi_{(i_1,\ldots, i_d)}.
	\]
Then $\Xi_j$ is a Hilbert--Schmidt operator if and only if $Q_j(2):=\sum_{i=1}^n |q_{ij}|^2<1$. Assume $Q_j(2)<1$ for each $j=1,\ldots,n$. Using the same identification as in Example 1, we can consider
	\[
		A=\sum_{j=1}^n \Xi_j\otimes E_{jj}\in L^2(M_n(W^*(X)\bar\otimes W^*(X)^{op})),
	\]
where $E_{jj}$ is the diagonal matrix with $jj$th entry 1 and all others zero. By \cite[Proposition 2]{NZ15b} $\J_X^*(A)=X$, and hence using \cite[Lemma 5.1]{NZ15} we see
	\[
		\Sigma^*(X\mid V_1)\leq \left(\sum_{j=1}^n \|\Xi_j - 1\otimes 1^{op}\|_{\tau\otimes\tau^{op}}^2\right)^{1/2} = \left(\sum_{j=1}^n \frac{Q_j(2)}{1-Q_j(2)}\right)^{1/2}.
	\]
\end{cex}


\begin{cex}{3: $q$-deformed Araki-Woods algebras \cite{Hia03}}

Once again we have $\mc{H}_\R=\R^n$, however we modify $\mc{H}$. Let $\{U_t\colon t\in\R\}$ be an orthogonal representation of $\R$ onto $\mc{H}_\R$, which extends to a one parameter group of unitary transformations on $\C^n$. By Stone's theorem, there exists a closed positive operator $B$ such that $U_t=B^{it}$ for every $t\in \R$. We then consider the following inner product
	\[
		\<f,g\>_U:=\<f, \frac{2}{1+B^{-1}}g\>,\qquad f,g\in\C^n,
	\]
and let $\mc{H}$ be the completion of $\C^n$ with respect to $\|\cdot\|_U$. Recall from \cite{Shl97} that this yields an isometric embedding of $\mc{H}_\R$ into $\mc{H}$ where the restriction of $\Re\<\cdot,\cdot\>_U$ to $\mc{H}_\R$ recovers the original inner product. In particular, our orthonormal basis $\{e_1,\ldots, e_n\}$ for $\mc{H}_\R$, are unit vectors inside in $\mc{H}$ with completely imaginary covariance.

Given $q\in[-1,1]$, one defines an inner product on $\mc{F}^{\text{finite}}$ (see \cite[Section 1]{Hia03}) so that for each $j=1,\ldots, n$
	\begin{align*}
		\ell(e_j)^*\Omega &=0\\
		\ell(e_j)^* f_1\otimes\cdots\otimes f_d &= \sum_{k=1}^d q^{k-1}\<f_k,e_j\>_U f_1\otimes \cdots \otimes \hat{f}_k\otimes\cdots\otimes f_d,\qquad f_1,\ldots, f_d\in \mc{H}.
	\end{align*}
The operators $x_j=\ell(e_j)+\ell(e_j)^*$ are called \emph{$q$-quasi-free semicircular operators}, and the vacuum state $\vphi$ (which is not a trace unless the orthogonal representation $\{U_t\colon t\in \R\}$ is trivial) is called a \emph{$q$-quasi free state}. Furthermore, the action of the modular operator $\Delta_\vphi$ is known:
	\[
		\Delta_\vphi x_j= \sum_{k=1}^n [B]_{jk} x_k,
	\]
where $[B]_{jk}$ is the $jk$th entry of the matrix representation of $B$ as a map on $\mc{H}_\R$ with respect to the basis $\{e_1,\ldots, e_n\}$. In the case $q=0$, the $x_j$ are Shlyakhtenko's \emph{quasi-free semicircular operators} and $\vphi$ is a \emph{free quasi-free state} \cite{Shl97}.

Letting $\Xi_q$ be exactly as in Example 1 (with $q^2n<1$ so it is a Hilbert--Schmidt operator), it was shown in \cite[Proposition 2.2]{Nel15} that
	\[
		\<x_j, p(X)\>_\vphi = \sum_{k=1}^n \< \left[\frac{2}{1+B}\right]_{jk} (\Delta_\vphi\otimes \Delta_{\vphi^{op}})(\Xi_q^*), [\partial_k p](X)\>_{\vphi\otimes\vphi^{op}},\qquad p\in\mathscr{P},
	\]
where $\left[\frac{2}{1+B}\right]_{jk}$ is the $jk$th entry of the matrix representation of $\frac{2}{1+B}$ as a map on $\mc{H}_\R$ with respect to the basis $\{e_1,\ldots, e_n\}$. Thus
	\[
		A=  (\Delta_\vphi\otimes \Delta_{\vphi^{op}})(\Xi_q^*)\otimes \frac{2}{1+B} \in L^2(M_n(W^*(X)\bar\otimes W^*(X)^{op}),(\vphi\otimes\vphi^{op})\circ\Tr)
	\]
is a free Stein kernel for $X$ relative to $V_1$. Using $\left[\frac{2}{1+B}\right]_{jj}=1$ we have
	\begin{align*}
		\Sigma^*(X\mid V_1)\leq \left(n\| (\Delta_\vphi\otimes \Delta_{\vphi^{op}})(\Xi_q^*)- 1\otimes 1^{op}\|_{\vphi\otimes\vphi^{op}}^2 + \sum_{j\neq k} \left|\left[\frac{2}{1+B}\right]_{jk}\right|^2 \| (\Delta_\vphi\otimes \Delta_{\vphi^{op}})(\Xi_q^*)\|_{\vphi\otimes\vphi^{op}}^2\right)^{1/2}.
	\end{align*}
It is not difficult to show that if $\{\xi_\ell^{(d)}\}_{\ell\in I_d}$ is an orthonormal basis for $\pi_d\mc{F}$, then
	\[
		\| (\Delta_\vphi\otimes \Delta_{\vphi^{op}})(\Xi_q^*)- 1\otimes 1^{op}\|_{\vphi\otimes\vphi^{op}}^2 = \sum_{d\geq 1}\sum_{\ell\in I_d} q^{2d} \|\Delta_\vphi^{-1} \xi_\ell^{(d)}\|_\vphi^2
	\]
Thus, in the estimate for $\Sigma^*(X\mid V_1)$ we see that the first term vanishes as $q\to 0$, whereas the sum vanishes as $B\to 1$. Since $W^*(X)$ is known to be a type $\mathrm{III}$ factor so long as $B\neq 1$, this would offer a counterexample to our conjecture were it not for the fact that
	\[
		A^\dagger=(\Delta_\vphi\otimes\Delta_{\vphi^{op}})(\Xi_q^*)\otimes \frac{2}{1+B^{-1}}\neq A.
	\]

\end{cex}


\section{One-dimensional case}
\label{section_1d}

In this section, we give the statements of our results in the one-dimensional case. The statements become simpler, in particular because we have an explicit expression for the entropy of a probability measure on $\R$: 
$$\chi(\mu)=\chi^*(\mu) := \int{\log|x - y|\mu(dx)\mu(dy)}+\frac{3}{4}+ \frac{1}{2}\log(2\pi)$$
while the relative entropy with respect to the semicircular law $d\sigma = \frac{1}{2\pi}\sqrt{4-x^2}dx$ (with unit variance) is given by
$$\chi(\mu| \sigma) := \frac{1}{2}\int{x^2 \mu(dx)} - \chi(\mu).$$
We note $\chi(\sigma)=\frac{1}{2}\log(2\pi e)$ and so $\chi(\sigma\mid \sigma)= -\frac{1}{2}\log(2\pi)$. The one-dimensional free Fisher information with respect to $\sigma$ is given by
$$\Phi^*(\mu|\sigma) = \int{(H\mu(x) - x)^2\mu(dx)}$$
where $H\mu(x) = \int{(x-y)^{-1}\mu(dy)}$ is the Hilbert transform of $\mu$. The LSI then takes the form
$$\chi(\mu|\sigma) + \frac{1}{2}\log(2\pi) \leq \frac{1}{2}\Phi^*(\mu|\sigma).$$
In addition to the general proof of Biane and Speicher \cite{BS01}, several other techniques have been successfully applied to obtain this one-dimensional inequality. Ledoux also proved it using a free analogue of the Brunn-Minkowski inequality \cite{Led05} (which itself was also proven using a random matrix approximation). Finally, Ledoux and Popescu \cite{LP09} gave a proof using mass transport, adapting arguments of Cordero-Erausquin \cite{CE02} in the classical case. 

The Stein factor $A$ associated to a probability measure (and relative to the semicircular law) is a function of two variables, such that
$$\int{xP(x)\mu(dx)} = \int{\int{A(x,y)\tilde{P}(x,y)\mu(dx)\mu(dy)}}$$
for all polynomials $P$, where $\tilde{P}$ is the free difference quotient, which takes value $\tilde{P}(x,y) = P'(x)$ if $x = y$, and $\frac{P(x) - P(y)}{x - y}$ otherwise. The measure $\mu$ is semicircular with unit variance iff $A$ is constantly equal to $1$, and the Stein discrepancy is given by
$$\Sigma^*(\mu\mid \sigma) := \left(\int{|A(x,y) - 1|^2\mu(dx)\mu(dy)}\right)^{1/2}.$$

In this context, the free HSI inequality takes the following form: 

\begin{cor}
For any centered measure $\mu$, given its Stein kernel $A$, we have
\begin{align*}
\frac{1}{2}\int{x^2\mu(dx)} &- \int{\log|x - y|\mu(dx)\mu(dy)} - \frac{3}{4}\\
&\leq \frac{1}{2}\int{|A(x,y) - 1|^2\mu(dx)\mu(dy)} \times \log\left(1 + \frac{\int{(H\mu(x) - x)^2\mu(dx)}}{\int{|A(x,y) - 1|^2\mu(dx)\mu(dy)}}\right).
\end{align*}
\end{cor}


\section{Rate of convergence in the free CLT}

As in the classical case, the free HSI inequality can be used to give bounds on how fast the entropy converges along the free Central Limit Theorem. Let us first recall the statement of the free CLT: 

\begin{thm}[Free CLT, \cites{Voi86, Sp90}]
For each $j=1,\ldots, N$ let $(x_j^{(N)})_{N\in \N}$ be a sequence of freely independent, identically distributed random variables. Assume each $x_j^{(N)}$ is centered and that the covariance of the $n$-tuples $X^{(N)}:=(x_1^{(N)},\ldots, x_n^{(N)})$ is the identity. If
$$Y^{(N)} := \frac{X^{(1)} + .. + X^{(N)}}{\sqrt{N}},$$
then $Y^{(N)}=(y_1^{(N)},\ldots, y_n^{(N)})$ converges (in moments) to a free $(0,1)$-semicircular $n$-tuple. 
\end{thm}

In \cite{CG13}, Chystiakov and Gotze showed that for a sequence of freely independent, identically distributed 1-tuples of random variables $(x_1^{(N)})_{N\in\N}$, the entropy of the renormalized convolution $y_1^{(N)}$ satisfies 
$$\chi^*(y_1^{(N)}|\sigma ) = \chi^*(\sigma|\sigma) + O(N^{-1})$$
under the assumptions that the random variables $x_1^{(N)}$ are bounded (actually, they show a more precise expansion, and identify the prefactor in the leading term of order $N^{-1}$ in the expansion). The order of magnitude of the reminder is sharp. The proof is based on an Edgeworth expansion of the distribution function of normalized sums. 

As in the classical case, we can use the HSI inequality to obtain a rate of convergence of order $\frac{\log N}{N}$ for random variables with finite free Fisher information and free Stein discrepancy. While the rate is slightly suboptimal, the technique has the advantage of avoiding technical computations, is easily applied to the multivariate case, and can be adapted to freely independent but not identically distributed $X_i$. 

A quantitative multi-dimensional free CLT was obtained by Speicher in \cite{Sp12}. The estimates obtained there involve distances between Cauchy transforms, and are obtained via a reduction to a one-dimensional operator-valued free central limit theorem. 

\begin{thm}
For each $j=1,\ldots, N$ let $(x_j^{(N)})_{N\in \N}$ be a sequence of freely independent, identically distributed random variables. Assume each $x_j^{(N)}$ is centered and that the covariance of the $n$-tuples $X^{(N)}:=(x_1^{(N)},\ldots, x_n^{(N)})$ is the identity. For $(a^{(N)}_1,\ldots, a_N^{(N)})_{N\in \N}$ an array of real numbers satisfying $\sum_{\ell=1}^N \hspace{1mm} (a_\ell^{(N)})^2 = 1$, let 
$$Y^{(N)} := a_1^{(N)}X^{(1)}+\cdots + a_N^{(N)}X^{(N)}.$$
Assume that the free Fisher information and the free Stein discrepancy of $X^{(1)}$ relative to $V_1$ are both finite. Let $\sigma_N := \sum_{\ell=1}^N \hspace{1mm} (a_\ell^{(N)})^4$. Then 
$$|\chi^*(Y^{(N)}|V_1) - \chi^*(S|V_1)| = O(\sigma_N \log (\sigma_N^{-1})).$$
\end{thm}
In particular, when $a_\ell^{(N)} = N^{-1/2}$ for each $\ell=1,\ldots, N$ we get a rate of order $\frac{\log N}{N}$ in the free CLT (instead of the sharp $\frac{1}{N}$ of \cite{CG13}). So in exchange for allowing non-identical weights in the sum and covering the multivariate case, we end up with a slightly worse rate. This situation is the same as for the classical entropic CLT obtained via the HSI inequality of \cite{LNP15}.

\begin{proof}
We first compute the free Stein kernel of $Y^{(N)}$ in terms of the kernels for $X^{(1)},\ldots,X^{(N)}$, which we denote $A_1,\ldots, A_N$. We may assume $\|A_\ell - 1\|_{HS}=\Sigma^*(X^{(\ell)}\mid V_1)=\Sigma^*(X^{(1)}\mid V_1)$ for each $\ell=1,\ldots, N$ (the latter equality following by the identical distributions of $X^{(1)},\ldots, X^{(N)}$). Since the $X^{(\ell)}$ are free, it follows that for $P\in \mathscr{P}^n$
	\begin{align*}
		\<Y^{(N)}, P(Y^{(N)})\>_2 &=\sum_{\ell=1}^N a_\ell^{(N)} \< X^{(\ell)}, P(Y^{(N)})\>_2\\
			&= \sum_{\ell=1}^N a_\ell^{(N)} \< A_\ell, a_\ell^{(N)} [\J P](Y^{(N)})\>_{HS}.
	\end{align*}
Thus if $\mc{E}_N$ is the projection onto $L^2(M_N(W^*(Y^{(N)})\bar\otimes W^*(Y^{(N)})^{op})$, then
	\[
		A:=\mc{E}_N\left(\sum_{\ell=1}^N \left(a_\ell^{(N)}\right)^2 A_\ell\right)
	\]
is a free Stein kernel for $Y^{(N)}$ relative to $V_1$. Therefore
	\begin{align*}
		\Sigma^*(Y^{(N)}\mid V_1)^2 &\leq \| A - 1\|_{HS}^2\\
			&= \left\| \mc{E}_N\left(\sum_{\ell=1}^N \left(a_\ell^{(N)}\right)^2(A_\ell - 1)\right)\right\|_{HS}^2\\
			&= \sum_{\ell}^N \left(a_\ell^{(N)}\right)^4 \|A_\ell -1\|_{HS}^2\\
			&= \sigma_N \Sigma^*(X^{(1)}\mid V_1)^2,
	\end{align*}
where the third line follows by freeness. Moreover, the Fisher information of $Y^{(N)}$ is smaller than that of $X^{(1)}$, as a direct consequence of the free Stam inequality (see \cite[Proposition 6.5]{VoiV}). The result then follows from a direct application of the HSI inequality to $Y^{(N)}$ and monotonicity of $t\mapsto t\log(1+c t^{-1})$ for $c,t\geq 0$.
\end{proof}


\section{Bound on the deficit in the log-Sobolev inequality}

We produce the free analogue of \cite[Theorem 1.1]{BGRS14}. Let $S=(s_1,\ldots, s_n)$ and $S'=(s'_1,\ldots, s'_n)$ each be free $(0,1)$-semicircular $n$-tuples in $M^n$, free from each other and free from $x_1,\ldots, x_n$. Fix $\rho>0$ and let $\frac{1}{\sqrt{\rho}}S=Y=(y_1,\ldots, y_n)$, a free $(0,\rho^{-1})$-semicircular $n$-tuple. For each $t\geq 0$ we denote
	\begin{align*}
		X(t)&=(x_1(t),\ldots, x_n(t))=(x_1+\sqrt{t}s'_1,\ldots, x_n+\sqrt{t}s'_n)\\
		Y(t)&=(y_1(t),\ldots, y_n(t))=(y_1+\sqrt{t}s'_1,\ldots, y_n+\sqrt{t}s'_n).
	\end{align*}
Note that $Y(t)$ is a free $(0,\rho^{-1}+t)$-semicircular $n$-tuple.

If $Z$ is a free $(0,\rho^{-1})$-semicircular $n$-tuple, then denote
	\[
		\Delta^*(X,Z)= \chi^*(X\mid V_\rho) - \chi^*(Z\mid V_\rho).
	\]
Also define $d(t)=t-\log(1+t)$ for $t>-1$.

\begin{thm}
	\[
		\Phi^*(X\mid V_\rho) - 2\rho \Delta^*(X,\sqrt{\rho}^{-1} S) \geq n\rho\cdot d\left( \frac{\Phi^*(X)}{n\rho} - 1\right).
	\]
\end{thm}
\begin{proof}
We compute using \cite[Proposition 7.5]{VoiV}
	\begin{align*}
		\Delta^*(X(t),Y(t))&= \chi^*(Y(t)) - \chi^*(X(t)) + \frac{1}{2(\rho^{-1}+t)}\sum_{j=1}^n \tau( x_j(t)^2 - y_j(t)^2)\\
			&= \chi^*(Y(t)) - \chi^*(X(t)) + \frac{1}{2(\rho^{-1}+t)}\sum_{j=1}^n \tau( x_j^2 - y_j^2)\\
			&= \chi^*(Y)-\chi^*(X) +\frac{1}{2}\int_0^t \Phi^*(Y(r)) -\Phi^*(X(r))\ dr + \frac{1}{2(\rho^{-1}+t)}\sum_{j=1}^n \tau( x_j^2 - y_j^2)\\
			&= \Delta^*(X,Y) + \frac{1}{2}\int_0^t \Phi^*(Y(r)) -\Phi^*(X(r))\ dr - \frac{\rho t}{2(\rho^{-1}+t)}\sum_{j=1}^n \tau( x_j^2 - y_j^2).
	\end{align*}
Equivalently,
	\[
		\Delta^*(X,Y) = \Delta^*(X(t), Y(t)) +\frac{1}{2}\int_0^r \Phi^*(X(r)) - \Phi^*(Y(r))\ dr + \frac{\rho t}{2(\rho^{-1}+t)} \sum_{j=1}^n \tau(x_j^2 - y_j^2).
	\]
From the free Stam inequality \cite[Proposition 6.5]{VoiV} we have
	\[
		\Phi^*(X(r)) \leq \frac{1}{\Phi^*(X)^{-1} + \Phi^*(\sqrt{r}S')^{-1}} = \frac{n\Phi^*(X)}{n+r\Phi^*(X)}.
	\]
Thus
	\begin{align*}
		\int_0^t \Phi^*(X(r)) - \Phi^*(Y(r))\ dr &\leq \int_0^r \frac{n\Phi^*(X)}{n+r\Phi^*(X)} - \frac{n}{\rho^{-1}+r}\ dr\\
			& = n \log\left( \frac{n+t\Phi^*(X)}{n(1+\rho t)}\right).
	\end{align*}
Continuing our previous computation we have
	\begin{align}\label{pre_limit_ineq}
		\Delta^*(X,Y) \leq  \Delta^*(X(t), Y(t)) + \frac{n}{2} \log\left( \frac{n+t\Phi^*(X)}{n(1+\rho t)}\right)+ \frac{\rho t}{2(\rho^{-1}+t)} \sum_{j=1}^n \tau(x_j^2 - y_j^2).
	\end{align}
We claim $\lim_{t\to\infty} \Delta^*(X(t),Y(t))=0$. By our previous computation, it suffices to show $\lim_{t\to\infty} \chi^*(Y(t)) - \chi^*(X(t))=0$. Note that $\chi^*(Y(t))= \frac{n}{2}\log(2\pi e (\rho^{-1}+t))$ and so by \cite[Proposition 7.2]{VoiV} and \cite[Proposition 7.5.(a)]{VoiV} we have
	\begin{align*}
		\frac{n}{2}\log\left(\frac{\rho^{-1}+t}{n^{-1}\Var(X) +t}\right) \leq \chi^*(Y(t)) - \chi^*(X(t)) \leq \frac{n}{2}\log\left(\frac{\rho^{-1}+t}{t}\right),
	\end{align*}
which establishes the claim.

Hence, taking the limit as $t\to\infty$ in (\ref{pre_limit_ineq}) yields
	\begin{align*}
		\Delta^*(X,Y) \leq \frac{n}{2}\log\left(\frac{\Phi^*(X)}{n\rho}\right) + \frac{\rho }{2} \left[\Var(X) - \frac{n}{\rho}\right].
	\end{align*}
Replacing $\rho\Var(X) - n $ with $\frac{1}{\rho}\left(\Phi^*(X\mid V_\rho) - \Phi^*(X)+n\rho\right)$ yields
	\begin{align*}
		\Delta^*(X,Y) &\leq \frac{n}{2}\log\left(  \frac{\Phi^*(X)}{n\rho}\right)+ \frac{1}{2\rho }\Phi^*(X\mid V_\rho) - \frac{1}{2\rho }\Phi^*(X) + \frac{n}{2}\\
			&= \frac{1}{2\rho} \Phi^*(X\mid V_\rho) - \frac{n}{2} d\left( \frac{\Phi^*(X)}{n\rho} - 1\right),
	\end{align*}
which is equivalent to the claimed inequality.
\end{proof}


\section{A Microstates Free log-Sobolev Inequality}

In this section, we prove a log-Sobolev inequality for the microstates free entropy. This result is simply a multi-variable analogue of \cite[Theorem 4]{LP09}. The proof of \cite[Theorem 4]{LP09} relies on the existence of transport between the measure of interest and the minimizer of the free entropy relative to a potential. Recall that by Lemma \ref{minimizer_is_free_Gibbs_state} the minimizer of the microstates free entropy relative to a potential $V\in \mathscr{P}^{(R)}$ has as a joint law a free Gibbs state with potential $V$. Not only will we have to assume the existence of free transport in the multi-variable case, which is already a very restrictive assumption, but moreover we must assume that the free transport map is monotone and implemented by an invertible $n$-tuple of power series.

\begin{defi}
For $f\in \mathscr{P}^{(R)}$, $R>0$, we say that $f$ is convex if $\J\D f \geq 0$ as an element of $M_n(\mathscr{P}\hat\otimes_R \mathscr{P}^{op})$.
\end{defi}

We have the following lemma, which is the analogue of the fact that a classical convex function lies above any of its tangent lines.

\begin{lem}\label{convex}
Let $(M,\tau)$ be a tracial von Neumann algebra and let $f\in\mathscr{P}^{(R)}$, $R>0$, be convex. If $A,B\in M^n$ (self-adjoint) satisfy $\|A\|,\|B\|\leq R$, then
	\[
		\tau(f(A)) \geq \tau(f(B)+ [\D f](B)\cdot (A-B)).
	\]
\end{lem}
\begin{proof}
The convexity of $f$ implies $[\J \D f](A,B)$ is positive in $M_n(M\bar\otimes M^{op})$. Hence
	\begin{align*}
		0&\leq \<[\J\D f](A,B)\#(A-B), (A-B)\>_2 = \sum_{j,k=1}^n \tau( (A_j-B_j)[\partial_k \D_j f](A,B)\# (A_k - B_k))\\
			&= \sum_{j=1}^n \tau( (A_j-B_j)\left\{ [\D_j f](A) - [\D_j f](B)\right\}) = \tau( [\D f](A)\cdot (A -B)- [\D f](B)\cdot (A-B)).
	\end{align*}
Now, for $t\in [0,1]$ define $A_t:=tA+(1-t)B$ and 
	\[
		h(t):=\tau( f(A_t) - f(B) - [\D f](B) \cdot (A_t - B)).
	\]
Note that $A_t - B = t(A-B)$. We clearly have $h(0)=0$ and we wish to show $h(1)\geq 0$, which we shall demonstrate by showing $h'(t)\geq 0$. We compute
	\begin{align*}
		\frac{d}{dt} h(t) &= \tau(\sum_{j=1}^n [\partial_j f](A_t)\#(A-B) - [\D f](B)\cdot (A-B))\\
			&= \frac{1}{t} \tau( [\D f](A_t)\cdot (A_t -B) - [\D f](B)\cdot (A_t -B)),
	\end{align*}
which is positive by our initial observation above.
\end{proof}

We also have a much weaker converse of Lemma \ref{minimizer_is_free_Gibbs_state}:

\begin{prop}
Let $Y=(y_1,\ldots, y_n)\in M^n$ have as a joint law a free Gibbs state with potential $V$ for some convex, self-adjoint $V\in \mathscr{P}^{(R)}$, $R> \|Y\|$. Let $F\in (\mathscr{P}^{(R)})^n$ be such that it has a compositional inverse $G\in (\mathscr{P}^{(\|F\|_R)})^n$ and $\J F\geq 0$. Then $\chi(F(Y)\mid V)\geq \chi(Y\mid V)$.
\end{prop}
\begin{proof}
Note that $Y=G(F(Y))$. We compute using Lemma \ref{convex} and \cite[Proposition 3.5]{VoiII} 
	\begin{align*}
		\chi(F(Y)\mid V) - \chi(Y\mid V) &= \tau(V(F(Y)) - V(Y)) - \chi(F(Y))+\chi(Y)\\
			&\geq \tau( [\D V](Y)\cdot (F(Y)- Y)) + (\tau\otimes\tau^{op})\circ\Tr( \log|[\J G](F(Y))|).
	\end{align*}
Now, since the joint law of $Y$ is a free Gibbs state with potential $V$, $[\D_1 V](Y),\ldots, [\D_n V](Y)$ are the conjugate variables to $Y$. So using (\ref{conj_var_equation}) we have
	\begin{align*}
			\chi(F(Y)\mid V) -\chi(Y\mid V) \geq (\tau\otimes \tau^{op})\circ\Tr( [\J F](Y) - 1 + \log| [\J G](F(Y))|).
	\end{align*}
Then, since $G\circ F$ is the identity, we have $[\J G](F(Y)) =[(\J F)^{-1}](Y)$. So continuing the above computation we have
	\begin{align*}
		\chi(F(Y)\mid V) - \chi(Y\mid V) \geq (\tau\otimes\tau^{op})\circ\Tr( [\J F](Y) - 1  - \log[\J F](Y)) \geq 0,
	\end{align*}
since $t- 1 -\log(t)\geq 0$ for $t\geq 0$.
\end{proof}

\begin{prop}
Suppose $Y=(y_1,\ldots, y_n)\in M^n$ has a free Gibbs state with potential $V$ as a joint law, where $V=V^*\in \mathscr{P}^{(R)}$, $R>\|Y\|$, is such that $V-\rho V_1$ is convex for some $\rho>0$. Suppose $F\in (\mathscr{P}^{(R)})^n$ satisfies:
	\begin{itemize}
		\item[(i)] there exists $G\in (\mathscr{P}^{(\|F\|_R)})^n$ a compositional inverse of $F$;
		\item[(ii)]  $\J G\geq 0$; and
		\item[(iii)] $F\in (\mathscr{P}^{(\|G\|_{\|F\|_R})})^n$.
	\end{itemize}
Then
	\[
		\chi(F(Y)\mid V) - \chi(Y\mid V) \leq \frac{1}{2\rho} \Phi^*(F(Y)\mid V)
	\]
\end{prop}
\begin{proof}
Denote $F(Y)=X=(x_1,\ldots, x_n)$. If $\Phi^*(X)=\infty$, then there is nothing to show. Assume $\Phi^*(X)<\infty$ and that $\xi_1,\ldots, \xi_n$ are the conjugate variables to $x_1,\ldots, x_n$. Denote $\Xi=(\xi_1,\ldots, \xi_n)$. Observe that the claimed inequality is equivalent to:
	\begin{align*}
		\frac{1}{2\rho} \| \Xi - [\D V](X) \|_2^2 &+ \tau( V(Y) - V(X) - [\D V](X)\cdot (Y-X))\\
									&- \<\Xi - [\D V](X), Y - X\>_2 + \< \Xi, Y- X\>_2 +\chi(X) - \chi(Y)\geq 0.
	\end{align*}
Let $W=V- \rho V_1$, so that by Lemma \ref{convex} we have $\tau( W(Y) - W(X) - [\D W](X)\cdot (Y-X)) \geq 0$ or
	\begin{align*}
		\tau( V(Y) - V(X) - [\D V](X) \cdot (Y - X)) &\geq \rho\tau( V_1(Y) -  V_1(X) -  X\cdot (Y - X))\\
			&=\frac{\rho}{2} \tau( Y\cdot Y - X\cdot X - 2X\cdot (Y- X))\\
			& = \frac{\rho}{2} \|Y-X\|_2^2
	\end{align*}
Noting that
	\[
		\frac{1}{2\rho} \|\Xi-[\D V](X)\|_2^2 + \frac{\rho}{2} \|Y-X\|_2^2 - \<\Xi- [\D V](X), Y-X\>_2 = \frac{1}{2\rho} \| \Xi - [\D V](X) - \rho(Y-X)\|_2^2,
	\]
we see that it suffices to show $\<\Xi, Y-X\>_2 + \chi(X) - \chi(Y)\geq 0$. Conditions (i)-(iii) imply that the full hypotheses of \cite[Proposition 3.5]{VoiII} are satisfied. Using this and (\ref{conj_var_equation}) (on $Y=G(X$)) we have
	\begin{align*}
		\<\Xi, Y-X\>_2 + \chi^*(X) - \chi^*(Y) = (\tau\otimes\tau^{op})\circ\Tr([\J G](X) - 1 -\log [\J G](X)),
	\end{align*}
which is positive since $t-1 - \log(t)\geq 0$ for $t\geq 0$.
\end{proof}


\bibliography{references}

\end{document}